\documentclass[10pt]{amsart}

\usepackage{pst-all}
\usepackage{amsmath}
\usepackage[latin1]{inputenc}
\usepackage{amsfonts}
\usepackage{amssymb}
\usepackage{multicol}
\usepackage{verbatim}
\usepackage{amsthm}
\usepackage{amscd}
\usepackage{pb-diagram}
\usepackage[mathscr]{eucal}
\usepackage[latin1]{inputenc}

\usepackage{amsmath}

\usepackage{amsfonts}

\usepackage{amssymb}

\usepackage{amsthm}

\usepackage{graphics}

\newcommand{\ZZ}{\mathbb{Z}}

\newcommand{\CC}{\mathbb{C}}

\newcommand{\NN}{\mathbb{N}}

\newcommand{\yl}{15pt}

\newcommand{\ffbox}[1]{
\setbox9=\hbox{$\scriptstyle\overline{1}$}
\framebox[\yl][c]{\rule{0mm}{\ht9}${\scriptstyle #1}$}
}

\newcommand{\RR}{\mathbb{R}}   

\newcommand{\Glie}{\mathfrak{g}}

\newcommand{\Hlie}{\mathfrak{h}}

\newcommand{\demo}{\noindent {\it \small Proof:}\quad}

\newcommand{\U}{\mathcal{U}}

\newtheorem{thm}{Theorem}[section]

\newtheorem{defi}[thm]{Definition}

\newtheorem{prop}[thm]{Proposition}

\newtheorem{rem}[thm]{Remark}

\usepackage[all]{xy}
\title{Quantum toroidal algebras and their representations}

\author{David Hernandez}
\address{CNRS - UMR 8100 : Laboratoire de Math\'ematiques de Versailles, 45 avenue des Etats-Unis , Bat. Fermat, 78035 VERSAILLES, 
FRANCE}
\email{hernandez @ math . cnrs . fr}
\urladdr{http://www.math.uvsq.fr/\textasciitilde hernandez}

\begin{document}

\begin{abstract} Quantum toroidal algebras (or double affine quantum algebras) are defined from quantum affine Kac-Moody algebras by using the Drinfeld quantum affinization process. They are quantum groups analogs of elliptic Cherednik algebras (elliptic double affine Hecke algebras) to whom they are related via Schur-Weyl duality. In this review paper, we give a glimpse on some aspects of their very rich representation theory in the context of general quantum affinizations. We illustrate with several examples. We also announce new results and explain possible further developments, in particular on finite dimensional representations at roots of unity.

\vskip 4.5mm

\noindent {\bf 2000 Mathematics Subject Classification :} Primary 17B37, Secondary 81R50, 82B23.

\noindent {\bf Key words :} Quantum toroidal algebras, integrable representations, elliptic Cherednik algebras.

\end{abstract}

\maketitle

\section{Introduction} Let us consider the following two important "generalizations" of finite dimensional complex simple Lie algebras $\Glie$ which have been intensively studied : Kac-Moody affine Lie algebras $\hat{\Glie}$ (infinite dimensional analogs of $\Glie$) and quantum groups $\U_q(\Glie)$ (quantization of the enveloping algebra $\U(\Glie)$). The combination of the two gives quantum affine algebras $\U_q(\hat{\Glie})$ which have a very rich representation theory. In fact $\U_q(\hat{\Glie})$ can be obtained by quantization of the enveloping algebra $\U(\hat{\Glie})$, or directly from $\U_q(\Glie)$ by the Drinfeld quantum affinization process :
$$\begin{CD}
\Glie  @>{\text{Affinization}}>>  \hat{\Glie}\\
@VV{\text{Quantization}}V     @VV{\text{Quantization}}V\\
\U_q(\Glie)  @>{\text{Drinfeld Quantum Affinization}}>>  \U_q(\hat{\Glie})
\end{CD}.$$
The Drinfeld Quantum Affinization can be applied to any quantum Kac-Moody algebra. From a quantum group we get a quantum affine algebra \cite{bec, Dri2}. In general we do not get a quantum Kac-Moody algebra, but a new class of algebras called quantum affinizations. For example from a quantum affine Kac-Moody algebra one gets a toroidal (or double affine) quantum algebra (it is not a quantum Kac-Moody algebra, and so it can not be affinized again by this process !). In type $A$ these algebras $\U_q(sl_{n+1}^{tor})$ were first introduced by Ginzburg-Kapranov-Vasserot \cite{gkv} and then in the general context \cite{jin, Naams}. In type $A$ they are in Schur-Weyl duality with elliptic Cherednik algebras \cite{vv1}, and so can be considered as quantum groups analogs of elliptic Cherednik algebras (elliptic double affine Hecke algebras). 

Their representation theory seems to be very rich and promising, and its study is just at the beginning. In this review we explain the construction of these algebras, and we explain some aspects of their representation theory. One of the first success of the theory is the action on the $q$-Fock space \cite{vv2, stu} of $\U_q(sl_{n+1}^{tor})$. In this review we will focus on the representation theory in the general context of quantum affinizations (in particular on the integrable representations as defined in \cite{Naams, m2, h2}). We explain that the Drinfeld "coproduct" gives a structure of tensor category (fusion tensor category) by using a renormalization process and review applications \cite{h5}. 

In the case of simply-laced Dynkin diagrams, a powerful geometric approach to the representation theory of quantum affinizations has been developed by Nakajima \cite{Naams} by using quiver varieties, and for the particular cases of simply-laced quantum toroidal algebras see \cite{npc, vv3}. As the picture is given in details in \cite[Section 6]{npc}, we refer to it. These constructions, as well as promising new directions \cite{Sc, Sc2, nag}, are only sketchily reviewed in the present article. For results on classical toroidal algebras, we refer to the recent review \cite{ra}, references therein and \cite{cl}. The relation between the two cases is studied in \cite[Section 4]{e} by using the results in \cite{kass, mry}.

We also announce in this paper new results, in particular we give an example from a new natural family of finite dimensional representations at roots of unity of quantum toroidal algebras, which are motivated by the monomial realization of Kashiwara extremal crystals. We also discuss a relation with quantum affinizations associated with infinite Dynkin diagrams which implies explicit character formulas for Kirillov-Reshetikhin modules of $\U_q(sl_{n+1}^{tor})$.

This review is organized as follows : in Section \ref{un} we review basic definitions on quantum toroidal algebras, we explain the Schur-Weyl duality with elliptic Cherednik algebras and give as an application the representation on the $q$-Fock space. In section \ref{deux} we review results which are valid for arbitrary quantum affinizations, in particular on categories of integrable representations. Section \ref{trois} includes new material, in particular on examples of a new class of representations and relations to "infinite" quantum affinizations. Other possible developments are also discussed.

The bibliography is quite long as we tried to give many relevant references for quantum toroidal algebras and related topics. We apologize for the papers which could be missing. The subject grew in different directions as several promising directions emerged and this is probably the last moment to make an attempt to describe all these different directions in a short journal paper. The author made several choices, for example by explaining basic examples and focusing in the second part on the results known for general quantum affinizations.

{\bf Acknowledgments :} This review was presented in a talk given at the conference on Cherednik algebras at Edinburgh in June 2007 and the author would like to thank the organizers for the invitation. Also he would like to thank B. Leclerc, M. Varagnolo, E. Vasserot and the referee for their comments, and K. Nagao for discussions on his results.

\section{Quantum toroidal algebras}\label{un}

In this section we give the general definition of quantum affinizations and in particular of quantum toroidal algebras. We also start to review some important properties.

\subsection{Generalized Cartan matrix} 

\begin{defi} A generalized Cartan matrix is a matrix $C=(C_{i,j})_{1\leq i,j\leq n}$ satisfying \label{carmat}
$C_{i,j}\in\mathbb{Z}$, $C_{i,i}=2$, $(i\neq j\Rightarrow C_{i,j}\leq 0)$ and $(C_{i,j}=0\Leftrightarrow C_{j,i}=0)$.
\end{defi}
We denote
$I=\{1,...,n\}$ and $l=\text{rank}(C)$. For $i,j\in I$, we set $\delta_{i,j}=0$ if $i\neq j$, and $\delta_{i,j}=1$ if $i=j$. In the following $C$ is fixed and supposed to be symmetrizable, that is to say there is a matrix
$D=\text{diag}(r_1,...,r_n)$ ($r_i\in\mathbb{N}^*$)\label{ri} such that $B=DC$\label{symcar} is symmetric. 

\noindent $q\in\mathbb{C}^*$ is not a root of unity and is fixed. We set $q_i=q^{r_i}$ and we set for $l\in\ZZ$ and $0\leq k\leq s$:
$$[l]_q=\frac{q^l-q^{-l}}{q-q^{-1}}\in\mathbb{Z}[q^{\pm}]\text{ , }[s]_q! = [s]_q\cdots [1]_q\text{ , }\begin{bmatrix}s\\k\end{bmatrix}_q = \frac{[s]_q!}{[s-k]_q![k]_q!}.$$
We fix a realization $(\mathfrak{h}, \Pi, \Pi^{\vee})$ of $C$ (see \cite{kac}): $\mathfrak{h}$ is a $2n-l$
dimensional $\mathbb{Q}$-vector space, $\Pi=\{\alpha_1,...,\alpha_n\}\subset \mathfrak{h}^*$ (set of the simple roots), $\Pi^{\vee}=\{\alpha_1^{\vee},...,\alpha_n^{\vee}\}\subset \mathfrak{h}$ (set of simple coroots) are defined such that for $1\leq i,j\leq n$, $\alpha_j(\alpha_i^{\vee})=C_{i,j}$.

\subsection{Quantum affinizations and quantum Kac-Moody algebras} 

The following algebra is associated to the data $q$ and $(\Hlie,\Pi,\Pi^\vee)$ :

\begin{defi}\label{defiqaf} The quantum affinization $\U_q(\hat{\Glie})$ is the algebra with generators $x_{i,m}^{\pm}$ ($i\in I, m\in\ZZ$), $k_h$ ($h\in \Hlie$), $h_{j,r}$ ($j\in I, r\in\ZZ-\{0\}$), central elements $c^{\pm 1/2}$ and relations :
$$k_hk_{h'}=k_{h+h'}\text{ , }k_0=1\text{ , }[k_{h},h_{j,m}] = 0,$$
\begin{equation}\label{h1}[h_{i,m},h_{j,m'}] = \delta_{m, - m'} \frac{[mB_{i,j}]_q}{m} \frac{c^m - c^{-m}}{q - q^{-1}},\end{equation}
$$k_{h}x_{j,r}^{\pm}k_{-h}=q^{\pm \alpha_j(h)}x_{j,r}^{\pm},$$
\begin{equation}\label{hd}[h_{i,m},x_{j,r}^{\pm}] = \pm \frac{1}{m}[mB_{i,j}]_q c ^{-|m|/2} x_{j,m+r}^{\pm},\end{equation}
$$[x_{i,r}^+,x_{j,r'}^-] = \delta_{i,j}\frac{c^{(r-r')/2}\phi^+_{i,r+r'}-c^{-(r-r')/2}\phi^-_{i,r+r'}}{q_i-q_i^{-1}},$$
\begin{equation}\label{hdd}x_{i,r+1}^{\pm}x_{j,r'}^{\pm} - q^{\pm B_{i,j}}x_{j,r'}^{\pm}x_{i,r+1}^{\pm}
=q^{\pm B_{i,j}}x_{i,r}^{\pm}x_{j,r'+1}^{\pm}-x_{j,r'+1}^{\pm}x_{i,r}^{\pm},\end{equation}
and for $i\neq j$, the Drinfeld-Serre relations :
$$\underset{\pi\in \Sigma_s}{\sum}\underset{k=0..s}{\sum}(-1)^k\begin{bmatrix}s\\k\end{bmatrix}_{q_i}x_{i,r_{\pi(1)}}^{\pm}...x_{i,r_{\pi(k)}}^{\pm}x_{j,r'}^{\pm}x_{i,r_{\pi(k+1)}}^{\pm}...x_{i,r_{\pi(s)}}^{\pm}=0,$$
where $s=1-C_{ij}$, $\Sigma_s$ is the symmetric group on $s$ letters, and :
$$\phi_i^\pm (z) = \underset{m\geq 0}{\sum}\phi_{i,\pm m}^{\pm}z^{\pm m} = k_{\pm r_i \alpha_i^{\vee}}\text{exp}(\pm(q-q^{-1})\underset{m'\geq 1}{\sum}h_{i,\pm m'}z^{\pm m'}).$$
\end{defi}

The subalgebra generated by the $k_h$, $x_i^{\pm} = x_{i,0}^\pm$ is denoted by $\U_q(\Glie)$. It is isomorphic to the quantum Kac-Moody algebra of generalized Cartan matrix $C$. That is why $\U_q(\hat{\Glie})$ is considered as a quantum affinization of $\U_q(\Glie)$. Note that it is well-known that $\U_q(\Glie)$ is a Hopf algebra with the coproduct :
$$\Delta(x_i^+) = x_i^+ \otimes 1 + k_i^{-1}\otimes x_i^+ \text{ , }\Delta(x_i^-) = x_i^- \otimes k_i + 1\otimes x_i^-\text{ , }\Delta(k_i) = k_i\otimes k_i.$$
For $J\subset I$, we denote by $\hat{\U}_J$ the subalgebra of $\U_q(\hat{\Glie})$ generated by the $x_{i,m}^{\pm}$ ($i\in J$, $m\in\ZZ$), the $k_{r_i\alpha_i^\vee}$ ($i\in J$) and the $h_{i,r}$ ($i\in J$, $r\in\ZZ^*$).

\noindent Examples : 

Suppose that $C$ is of finite type : 
$$\forall 1\leq m\leq n, \text{det}((C_{i,j})_{1\leq i,j\leq m}) > 0.$$
Then $\U_q(\Glie)$ is a quantum group, $\U_q(\hat{\Glie})$ is an untwisted quantum affine algebra \cite{bec, Dri2}.

Suppose that $C$ is of affine type : 
$$\forall 1\leq m\leq n - 1, \text{det}((C_{i,j})_{1\leq i,j\leq m}) > 0\text{ and }\text{det}(C) = 0.$$
Then $\U_q(\Glie)$ is a quantum affine algebra, $\U_q(\hat{\Glie})$ is a quantum toroidal algebra (double affine quantum algebra). $\U_q(\Glie) = \U_q^h(\Glie)\subset\U_q(\hat{\Glie})$ is called the horizontal quantum affine subalgebra of $\U_q(\hat{\Glie})$. The vertical quantum affine subalgebra is $\U_q^v(\Glie)$ generated by the $x_{i,m}^{\pm}, k_i^{\pm 1},h_{i,r}$ for $i$ not equal to the additional node of $\Glie$ (see \cite{kac}). When $\Glie$ is untwisted $\U_q^v(\Glie)$ is isomorphic to $\U_q(\Glie)$. Note that a quantum toroidal algebra is not isomorphic to a quantum Kac-Moody algebra, and so can not be affinized again by this process.

\begin{rem}\label{tp} In the particular case $\Glie = \hat{sl}_{n+1}$, an additional parameter $d\in\CC^*$ is added in \cite{gkv} (two parameter deformation). It fits well into the Schur-Weyl duality discussed bellow : in addition to $C$ we also consider a matrix $M$ :
$$C = \begin{pmatrix}2&-1&&0&-1\cr
-1&2&\cdots&0&0\cr
&\vdots&\ddots&\vdots&\cr
0&0&\cdots&2&-1\cr
-1&0&&-1&2
\end{pmatrix} \text{ , }M = \begin{pmatrix}
0&-1&&0&1\cr
1&0&\cdots&0&0\cr
&\vdots&\ddots&\vdots&\cr
0&0&\cdots&0&-1\cr
-1&0&&1&0
\end{pmatrix}.$$
Then $\U_{q,d}(\hat{\Glie}) = \U_{q,d}(sl_{n+1}^{tor})$ is the algebra with the same generators and relations as $\U_q(\hat{\Glie})$ except the following $d$-deformed relations instead of respectively relations (\ref{h1}), (\ref{hd}), (\ref{hdd}) (the other relations are not modified) :
$$[h_{i,m},h_{j,m'}] d^{m M_{i,j}} = \delta_{m, - m'} \frac{[mB_{i,j}]_q}{m} \frac{c^m - c^{-m}}{q - q^{-1}}\text{ for $m > 0$, $m'\neq 0$,}$$
$$[h_{i,m},x_{j,r}^{\pm}]d^{mM_{i,j}}=\pm \frac{1}{m}[mB_{i,j}]_q c ^{-|m|/2} x_{j,m+r}^{\pm},$$
$$d^{M_{i,j}}x_{i,r+1}^{\pm}x_{j,r'}^{\pm}-d^{M_{i,j}}q^{\pm B_{ij}}x_{j,r'}^{\pm}x_{i,r+1}^{\pm}
=q^{\pm B_{ij}}x_{i,r}^{\pm}x_{j,r'+1}^{\pm}-x_{j,r'+1}^{\pm}x_{i,r}^{\pm}.$$
The horizontal subalgebra generated by the $k_h$, $x_i^\pm$ is isomorphic $\U_q(\hat{sl}_{n+1})$ (it does not depend on $d$). The vertical subalgebra generated by the $d^{-im}x_{i,m}^{\pm}$, $k_i^{\pm 1}$, $d^{-ir}h_{i,r}$ for $1\leq i\leq n$ is isomorphic to $\U_q(\hat{sl}_{n+1})$ (for this point it is crucial that $\Glie = \hat{sl}_{n+1}$, and it is not known to the author how to define for general type such a $q,d$-deformation so that we have a vertical subalgebra isomorphic to the quantum affine algebra).\end{rem}

\subsection{Weight lattice and finite dimensional representations}\label{wlfd}

Let us denote by $\Lambda_1,...,\Lambda_n\in\mathfrak{h}^*$ (resp. $\Lambda_1^{\vee},...,\Lambda_n^{\vee}\in\mathfrak{h}$) the fundamental weights (resp. coweights). By definition $\alpha_i(\Lambda_j^{\vee})=\Lambda_i(\alpha_j^{\vee})=\delta_{i,j}$. Denote $P=\{\lambda \in\mathfrak{h}^*|\forall i\in I, \lambda(\alpha_i^{\vee})\in\mathbb{Z}\}$ the set of integral weights and $P^+=\{\lambda \in P|\forall i\in I, \lambda(\alpha_i^{\vee})\geq 0\}$ the set of dominant weights. For example we have $\alpha_1,...,\alpha_n\in P$ and $\Lambda_1,...,\Lambda_n\in P^+$. Denote $Q={\bigoplus}_{i\in I}\mathbb{Z} \alpha_i\subset P$ the root lattice and $Q^+={\sum}_{i\in I}\mathbb{N} \alpha_i\subset Q$. For $\lambda,\mu\in \mathfrak{h}^*$, write $\lambda \geq \mu$ if $\lambda-\mu\in Q^+$.

For a quantum group $\U_q(\Glie)$ the (type $1$) simple finite dimensional representations are parameterized by $P^+$ \cite{ro, lu}. For $\lambda\in P^+$, $V(\lambda)$ is the corresponding simple representation. For a quantum affine algebra $\U_q(\hat{\Glie})$, the (type $1$) simple finite dimensional representations are parameterized by $n$-tuples of (Drinfeld) polynomials $\mathcal{P}=(P_1(z),\cdots, P_n(z))$ where $P_i(z)\in\CC[z]$ and $P_i(0) = 1$ \cite{cp1}. For such a $\mathcal{P}$, $V(\mathcal{P})$ is the corresponding simple representation. These results will be stated in a more precise way latter in the general context of quantum affinizations.

\subsection{Schur-Weyl duality and Cherednik algebras} We suppose that $C$ is of affine type $A_n^{(1)}$ and that $c = 1$ in this section. We have the quantum group $\U_q(sl_{n+1})$, the quantum affine algebra $\U_q(\hat{sl}_{n+1})$, and the quantum toroidal algebra $\U_{q,d}(sl_{n+1}^{tor})$ as defined in Remark \ref{tp}. Let $x\in\CC^*$ and $\epsilon = q^2$.

\begin{defi} The double affine Hecke algebra (DAHA) $\mathbb{H}_l(\epsilon,x)$ is the algebra with generators $\sigma_i^{\pm 1}$, $X_j^{\pm 1}$, $Y_j^{\pm 1}$ ($1\leq i\leq l-1$, $1\leq j\leq l$) and relations :
$$(\sigma_i + q^{-1})(\sigma_i - q) = 0\text{ , }\sigma_i\sigma_{i+1}\sigma_i = \sigma_{i+1}\sigma_i \sigma_{i+1},$$ 
$$\sigma_i\sigma_j = \sigma_j \sigma_i\text{ if $|i - j| > 1$},$$
$$X_0Y_1 = x Y_1X_0\text{ , }X_2Y_1^{-1}X_2^{-1}Y_1 = \sigma_1^2,$$    
$$X_iX_j = X_j X_i\text{ , }Y_i Y_j = Y_j Y_i\text{ for $1\leq i,j\leq l$},$$
$$X_j \sigma_i = \sigma_i X_j\text{ , }Y_j \sigma_i = \sigma_i Y_j\text{ if $j\neq i, i+1$},$$
$$\sigma_iX_i\sigma_i =  X_{i+1}\text{ , }Y_i = \sigma_iY_{i+1}\sigma_i\text{ for $1\leq i\leq l-1$,}$$   
where $X_0 = X_1\cdots X_l$.
\end{defi}

\noindent These algebras were introduced by Cherednik and have remarkable applications \cite{che2, che3}. DAHA (defined here in their elliptic form) and related algebras (trigonometric, rational forms) are currently intensively studied (see for example the papers in the special volume of the conference on Cherednik algebras).

\noindent The subalgebra generated by the $\sigma_i$ is isomorphic to the Hecke algebra $\mathcal{H}_l(\epsilon)\simeq \mathcal{H}_l(\epsilon^{-1})$. The subalgebra generated by the $\sigma_i^{\pm 1}, X_i^{\pm 1}$ (resp. the $\sigma_i^{\pm 1}, Y_i^{\pm 1}$) is isomorphic to the affine Hecke algebra $\tilde{\mathcal{H}}_l(\epsilon)\simeq \tilde{\mathcal{H}}_l(\epsilon^{-1})$.
  
Consider $V_1$ the $(n+1)$-dimensional representation of $\U_q(sl_{n+1})$ corresponding to the fundamental weight $\Lambda_1$ : $V_1 = \bigoplus_{i= 0 \cdots n} \CC v_i$ where ($v_{-1} = v_{n+1} = 0$) :
$$K_j.v_i = q^{\delta_{j-1,i} - \delta_{j,i}} v_i\text{ , }x_j^+.v_i = \delta_{i,j}v_{i - 1}\text{ , }x_j^-.v_i = \delta_{i+1,j} v_{i + 1}.$$
A completely reducible representation of $\U_q(sl_{n+1})$ is said to be of level $l$ if each of its irreducible components is isomorphic to some irreducible component of $(V_1)^{\otimes l}$. Note that any finite dimensional representation of $\U_q(sl_{n+1})$ is completely reducible. A simple finite dimensional representation $V(\sum_{i = 1\cdots n}\lambda_i \Lambda_i)$ is of level $l$ if and only if $\sum_{i=1\cdots n}i \lambda_i = l$.

\noindent $(V_1)^{\otimes l}$ has a natural structure of $\mathcal{H}_l(\epsilon)$ left-module, and for $M$ a right $\mathcal{H}_l(q^2)$-module, $I(M) = M\otimes_{\mathcal{H}_l(\epsilon)} (V_1)^{\otimes l}$ has a natural structure of $\U_q(sl_{n+1})$-module. 

\begin{thm}\cite{j} $I$ is a functor from the category of finite dimensional right $\mathcal{H}_l(\epsilon)$-modules to the category of finite dimensional $\U_q(sl_{n+1})$-modules of level $l$. Moreover if $l\leq n$, $I$ is an equivalence of category.
\end{thm}

\begin{thm}\label{cp2t}\cite{cp2} There is a functor $\mathcal{F}$ from the category of finite dimensional right $\tilde{\mathcal{H}}_l(\epsilon)$-modules to the category of finite dimensional $\U_q(\hat{sl}_{n+1})$-modules of level $l$, which coincides with $I$ for the $\U_q(sl_{n+1})$-module structure. Moreover if $l\leq n$, $\mathcal{F}$ is an equivalence of category.
\end{thm}

\begin{thm}\cite{vv1} For a right $\mathbb{H}_l(\epsilon,(d^{-1}q)^{n+1})$-module $M$, the action of the algebra $\U_q(\hat{sl}_{n+1})$ on $\mathcal{F}(M)$ can be extended to an action of $\U_{q,d}(sl_{n+1}^{tor})$.
\end{thm}

\noindent A reversed functor is also defined in \cite{vv1}. Recently analog results were studied in \cite{g, g2} in the trigonometric and rational cases, where affine Yangians and new algebras (deformed double current algebras) are used in the duality and have the PBW-property. 

Let us look at the duality more precisely in the affine case \cite[Section 12.3]{cp1} : consider the affine Hecke algebra $\tilde{\mathcal{H}}_l(\epsilon)$ with generators $\sigma_i, z_i$. In the following all considered $\tilde{\mathcal{H}}_l(\epsilon)$-modules are right-modules.

\noindent For $A=(a_1,\cdots,a_l)\in (\CC^*)^l$ let $M_A=(\sum_j (z_j-a_j.1)\tilde{\mathcal{H}}_l(\epsilon)
)\setminus \tilde{\mathcal{H}}_l(\epsilon)$. Every finite dimensional $\tilde{\mathcal{H}}_l(\epsilon)$-module is a quotient of some $M_A$. The natural map $\mathcal{H}_l(\epsilon)\rightarrow M_A$ is an isomorphism of $\mathcal{H}_l(\epsilon)$-module and in particular $\text{dim}(M_A) = l!$. $M_A$ is reducible if and only if $a_j=\epsilon a_k$ for some $j,k$.
 
\noindent The $q$-segment of center $a$ and length $l$ is defined as $\{aq^{1-l},aq^{3-l},\cdots,aq^{l-1}\}$. Let $S=\{s_1,\cdots, s_p\}$ be a collection of $q$-segments whose sum of length is $l$ and elements are given by $A$. By using an element $C_S\in \mathcal{H}_l(\epsilon)$ explicitly defined with Kazhdan-Lusztig polynomials (see \cite{h}), one can consider $I_S=C_S.\mathcal{H}_l(\epsilon)\subset M_A$. Then $I_S$ is a $\tilde{\mathcal{H}}_l(\epsilon)$-submodule of $M_A$ with a unique irreducible subquotient $V_S$ in which $C_S$ has non zero image. Every finite dimensional simple $\tilde{\mathcal{H}}_l(\epsilon)$ module is isomorphic to one $V_S$.
 
\begin{thm}\cite{cp2} Let $l\leq n$. For $S=\{s_1,\cdots,s_r\}$ a collection of $q$-segments, let $a_r$ be the center and $l_r$ the length of $s_r$. Then $\mathcal{F}(V_S)\simeq V(\mathcal{P})$ where $\mathcal{P}=(P_i)_{i\in I}$ is defined by : 
$$P_i(u)=\prod_{\{r|l_r=i\}}(ua_r-1).$$
\end{thm}
\noindent (The $C_S$ are not given explicitly above but this Theorem gives a characterization of the representation $V_S$).
 
$\tilde{\mathcal{H}}_{l_1+l_2}(\epsilon)$ is a $\tilde{\mathcal{H}}_{l_1}(\epsilon)\otimes\tilde{\mathcal{H}}_{l_2}(\epsilon)$-module thanks to the unique algebra morphism 
$$i_{l_1,l_2} : \tilde{\mathcal{H}}_{l_1}(\epsilon)\otimes\tilde{\mathcal{H}}_{l_2}(\epsilon)\rightarrow\tilde{\mathcal{H}}_{l_1+l_2}(\epsilon),$$ 
$$i_{l_1,l_2}(\sigma_i\otimes 1) = \sigma_i\text{ , }i_{l_1,l_2}(z_j\otimes 1) = z_j\text{ , }i_{l_1,l_2}(1\otimes\sigma_i)=\sigma_{i+l_1}\text{ , }i_{l_1,l_2}(1\otimes z_j)=z_{j+l_1}.$$
For $M_1$ (resp. $M_2$) a $\tilde{\mathcal{H}}_{l_1}(\epsilon)$ (resp. $\tilde{\mathcal{H}}_{l_2}(\epsilon)$) module, the $\tilde{\mathcal{H}}_{l_1+l_2}(\epsilon)$-module $M_1\otimes_Z M_2$ (Zelevinsky tensor product) is defined by :
$$M_1\otimes_Z M_2 = (M_1\otimes M_2)\bigotimes_{\tilde{\mathcal{H}}_{l_1}(\epsilon)\otimes\tilde{\mathcal{H}}_{l_2}(\epsilon)}\tilde{\mathcal{H}}_{l_1+l_2}(\epsilon).$$
We have the following compatibility property with the usual tensor product for quantum affine algebras :
\begin{prop}\label{compt}\cite{cp2} $\mathcal{F}(M_1\otimes_Z M_2)\simeq \mathcal{F}(M_1)\otimes_\CC\mathcal{F}(M_2)$.\end{prop}

Remark : for $a\in\CC^*$, $M_{(a)} = V_{(a)}$ is of dimension $1$, and for  $(a_1,\cdots,a_l)\in(\CC^*)^l$, we have : $$M_{(a_1,\cdots, a_l)}\simeq V_{(a_1)}\otimes_Z\cdots \otimes_Z V_{(a_l)}.$$

\subsection{Representations of quantum toroidal algebras and the $q$-Fock space} One of the main achievement in the representation theory of quantum toroidal algebras (see \cite{vv2, stu}) is the existence of an action on the $q$-Fock space defined in the following way : 
$\Lambda$ is the vector space of basis indexed by $|j_1,j_2,\cdots >\in\ZZ^\NN$ such that 
$$j_1 > j_2 > \cdots$$
$$j_{k+1} = j_k - 1\text{ for $k>>0$.}$$
For $M\in\ZZ$, let $\Lambda^{(M)}$ be the subspace of $\Lambda$ generated by the $|j_1,j_2,\cdots >$ such that $j_k = M-k +1$ for $k>>1$. It is remarkable that the space $\Lambda^{(M)}$ has two structures of representation of the quantum affine algebra $\U_q(\hat{sl}_{n+1})$, one of level $0$ and one of level $1$ (see \cite{tuu, jkkmp} and references therein). Then these two representations can be glued together as a representation of the whole quantum toroidal algebra (see \cite{vv2, stu} for the detailed definition of the action, see also \cite{tu} for other levels and \cite{tuu, tu} for complements). The two actions of the quantum affine algebra correspond to the action respectively of the vertical and horizontal quantum affine subalgebras. By lack of space we do not give the precise definition of the representation as it is already described in several papers mentioned above.

\noindent One remarkable point in the construction explained in \cite{vv2} is that the Dunkl-Cherednik representation of $\mathbb{H}_l(q^2,x)$ on $\CC[z_1^{\pm 1},\cdots,z_m^{\pm 1}]$ (see \cite{che2, che3}) plays a crucial role as it is used through Schur-Weyl duality as a building piece of the $q$-Fock space representation. 

For vertex algebras constructions of representations of quantum toroidal algebras see \cite{fj, gj, gj2, Sa}. Additional results on the structure of quantum toroidal algebras are given in \cite{m1, m4, m5}. General results valid for representations of more general quantum affinizations will be reviewed in the following.

\subsection{Geometric approaches} 

Beyond the case of $ADE$ quantum affine algebras, the powerful geometric construction of representations via quiver varieties can be extended to simply-laced quantum affinizations \cite{Naams}. For the particular cases of simply-laced quantum toroidal algebras see \cite{npc, vv3}. We refer to the review paper \cite{npc}. One important point is to replace the top degree homology group in the study of a quiver variety corresponding to an affine Dynkin diagram by equivariant K-groups. As explained in \cite{Naams, npc}, the geometric approach also includes the root of unity case. In fact quiver varieties are also known to have direct relations to Cherednik algebras (see \cite{go}).

More recently in \cite{nag}, an isomorphism between the geometric construction and the construction via Schur-Weyl duality of the level-(0, 1) $q$-Fock space representation of the ($d$-deformed) quantum toroidal algebra is given. In particular simultaneous eigenvectors on the $q$-Fock space are constructed by using non symmetric Macdonald polynomials and the corresponding eigenvalues are computed.

One of the most interesting feature of the representation theory of quantum Kac-Moody algebra is the existence of canonical/crystal bases. In \cite{Sc, Sc2} the construction of canonical bases for quantum affinization of quantum Kac-Moody algebras is proposed (in the case of stared Dynkin diagrams). Categories of coherent sheaves on a smooth projective curve $X$, equivariant with respect to a fixed subgroup $G$ of $Aut(X)$ such that $X/G\simeq \mathbb{P}_1$, are used. This implies the existence of a canonical basis in quantum toroidal algebras of type $D_4^{(1)}$, $E_6^{(1)}$, $E_7^{(1)}$, $E_8^{(1)}$.
                      
\section{Integrable representations and category $\mathcal{O}$}\label{deux}

In this section we will focus on results for general quantum affinizations $\U_q(\hat{\Glie})\supset \U_q(\Glie)$ with trivial central charge $c^{\pm 1/2} = 1$.

\subsection{Integrable representations} For $V$ a representation of $\U_q(\Glie)$ and $\nu\in P$, the weight space $V_{\nu}$ of $V$ is :
$$V_\nu = \{v\in V|k_h.v = q^{\nu(h)}v,\forall h\in \Hlie\}.$$

\begin{defi} $V$ is said to be in the category $\mathcal{O}$ if 

i) $V = \bigoplus_{\nu\in P} V_\nu$,

ii) $V_\nu$ is finite dimensional for all $\nu$,

iii) $\{\nu|V_\nu\neq\{0\}\}\subset \bigcup_{j=1\cdots N}\{\nu|\nu\leq \lambda_j\}$ for some $\lambda_1,\cdots,\lambda_N\in P$.
\end{defi}

\noindent For example for $\lambda\in P$, a representation $V$ is said to be of highest weight $\lambda$ if there is $v\in V_\lambda$ such that $\forall i\in I, x_i^+.v = 0$ and $\U_q(\Glie).v = V$. Such a representation is in the category $\mathcal{O}$. For each $\lambda\in P$ there is a unique simple highest weight module $V(\lambda)$ of highest weight $\lambda$.

\begin{defi} $V$ is said to be integrable if :

i) $V = \bigoplus_{\nu\in P} V_\nu$,

ii) $V_\nu$ is finite dimensional for all $\nu$,

iii) $V_{\nu \pm N\alpha_i} = \{0\}$ for all $\nu\in P$, $N>>0$, $i\in I$.
\end{defi}

\noindent Then we have :

\begin{thm}\cite{lu} 
The simple integrable representations of $\U_q(\Glie)$ in the category $\mathcal{O}$ are the $(V(\lambda))_{\lambda\in P^+}$.
\end{thm}

\noindent In the case of quantum groups of finite type we get the finite dimensional representations discussed in Section \ref{wlfd}.

\begin{defi} A representation of $\U_q(\hat{\Glie})$ is said to be integrable (resp. in the category $\mathcal{O}$) if is integrable (resp. in the category $\mathcal{O}$) as a $\U_q(\Glie)$-module.\end{defi}
      
\noindent Remark : for quantum affine algebras, we get two different notion of integrability (as a quantum affine algebra can be seen as a quantum Kac-Moody algebra or as a quantum affinization).
      
We have a triangular decomposition of $\U_q(\hat{\Glie})$. For a quantum affine algebra it is proved in \cite{bec}, and in full generality in \cite{h2} (the crucial point is to prove the compatibility with the Drinfeld-Serre relations):

\begin{thm}\label{dtrian} We have an isomorphism of vector spaces :
$$\U_q(\hat{\Glie})\simeq \U_q(\hat{\Glie})^-\otimes\U_q(\hat{\Hlie})\otimes\U_q(\hat{\Glie})^-,$$ 
where $\U_q(\hat{\Glie})^{\pm}$ (resp. $\U_q(\hat{\Hlie})$) is generated by the $x_{i,m}^{\pm}$ (resp. the $k_h$, the $h_{i,r}$).
\end{thm}

\begin{defi} A representation $V$ of $\U_q(\hat{\Glie})$ is said to be of $l$-highest weight if there is $v\in V$ such that 

i) $V = \U_q(\hat{\Glie})^-.v$, 

ii) $\U_q(\hat{\Hlie}).v=\CC v$, 

ii) for any $i\in I, m\in\ZZ$, $x_{i,m}^+.v=0$.
\end{defi}

For $\gamma\in \U_q(\hat{\Hlie})\rightarrow \CC$ an algebra morphism, by Theorem \ref{dtrian} we have a corresponding Verma module $M(\gamma)$ and a simple representation $L(\gamma)$. For a representation $V$, we also denote by $V_\gamma$ the common Jordan block of elements of $\U_q(\hat{\Hlie})$ with eigenvalues given by $\gamma$ (it is called a pseudo $l$-weight space). The following result was proved for $\Glie$ of finite type in \cite{cp1}, of type $A_n^{(1)}$ in \cite{m2}, of general simply-laced type in \cite{Naams}. In the general case, the proof \cite{h2} generalizes the proof of \cite{cp1} :

\begin{thm}\label{cond} The simple integrable representations of $\U_q(\hat{\Glie})$ in the category $\mathcal{O}$ are the $L(\gamma)$ such that there is $(\lambda,(P_i)_{1\leq i\leq n})\in P\times (1+u\mathbb{C}[u])^{\times n}$ satisfying $\gamma(k_h) = q^{\lambda(h)}$ and for $i\in I$ the relation in $\mathbb{C}[[z]]$ (resp. in $\mathbb{C}[[z^{-1}]]$):
$$\gamma(\phi_i^\pm(z))=q_i^{\pm \text{deg}(P_i)}\frac{P_i(zq_i^{-1})}{P_i(zq_i)}.$$
\end{thm}
\noindent By analogy to quantum affine algebras, the $P_i$ are called Drinfeld polynomials. Usually $\lambda$ is not used in the case of quantum affine algebras as $\text{dim}(\Hlie) = n$ and so $\lambda$ is uniquely determined by the $\gamma(k_{r_i\alpha_i^\vee})$, that is to say by the $\text{deg}(P_i)$.
    
Examples : 

1) for $k\geq 0$, $a\in\CC^*$, $1\leq i\leq n$, the Kirillov-Reshetikhin module $W_{k,a}^{(i)}$ is the simple representation of weight $\Lambda_i$ with the $n$-tuple 
\begin{equation*}
P_j(u) = 
\begin{cases}
(1-ua)(1-uaq_i^2)\cdots(1-uaq_i^{2(k - 1)})\text{ for $j = i$,}
\\1\text{ for $j\neq i$.} 
\end{cases}
\end{equation*}

2) the $V_i(a) = W_{1,a}^{(i)}$ are called fundamental representations.

3) for $\lambda\in P$ satisfying ($\forall j\in I, \lambda(\alpha_j^{\vee})=0$), $L_\lambda$ is the simple module of weight $\lambda$ with the $n$-uplet satisfying $P_j(u)=1$ for all $j\in I$. It is of dimension $1$.
      
\subsection{Fusion tensor category} $\U_q(\Glie)$ is a Hopf algebra, but $\U_q(\hat{\Glie})$ is not a Hopf algebra in general. However one can combinatorially define a natural ring structure \cite{h2} (fusion product) on the Grothendieck group $\text{Rep}(\U_q(\hat{\Glie}))$ of integrable representations in the category $\mathcal{O}$. Indeed consider \cite{h2} a generalization of the Frenkel-Reshetikhin $q$-character morphism \cite{Fre} (see also \cite{kn, cm, e}) : for $V$ a representation in $\text{Rep}(\U_q(\hat{\Glie}))$ we put :
$$\chi_q(V) = \sum_\gamma \text{dim}(V_\gamma) e(\gamma),$$
where the sum is over all algebra morphism $\gamma : \U_q(\hat{\Glie})\rightarrow\CC$ and $e(\gamma)$ is a formal element. We get an injective group morphism :
$$\chi_q :\text{Rep}(\U_q(\hat{\Glie}))\rightarrow \mathcal{E}$$ 
where $\mathcal{E}$ is a completion of $\bigoplus_{\gamma: \U_q(\hat{\Hlie})\rightarrow \CC}\CC e(\gamma)$ and 
$\text{Rep}(\U_q(\hat{\Glie}))$ is the Grothendieck group of integrable representations in the category $\mathcal{O}$. 

In the rest of this section we suppose that the quantized Cartan matrix 
$$C_{i,j}(z) = [C_{i,j}]_z + \delta_{i,j}([2]_{z_i} - [2]_z)$$ 
is invertible (here we denote $z_i = z^{r_i}$). It is proved in \cite{h1} that 
$$\forall i,j\in I, (C_{i,j} < -1 \Rightarrow -C_{j,i} \leq r_i)$$ 
is a sufficient condition. This is satisfied for quantum toroidal algebras (in the case $A_1^{(1)}$ we have to use the convention $r_0 = r_1 = 2$). In this situation there is a characterization of $\text{Im}(\chi_q)$ \cite{h2} which generalizes the characterization in the case of quantum affine algebras \cite{Fre, Fre2}. Moreover we have a natural product
$$(\gamma \gamma')(k_h) = \gamma(k_h)\gamma(k_h)\text{ , }(\gamma\gamma')(\phi_i^\pm(z)) = \gamma(\phi_i^\pm(z))\gamma'(\phi_i^\pm(z)),$$
which induces a ring structure on $\mathcal{E}$. In the case of quantum affine algebras, there is a natural ring structure on $\text{Rep}(\U_q(\hat{\Glie}))$ induced by the tensor product and $\chi_q$ is a ring morphism \cite{Fre}. In general we do not have {\it a priori} a ring structure on $\text{Rep}(\U_q(\hat{\Glie}))$, but :

\begin{thm}\cite{h2} $\chi_q(\text{Rep}(\U_q(\hat{\Glie})))$ is a subring of $\mathcal{E}$. The product $*$ induced on $\text{Rep}(\U_q(\hat{\Glie}))$ by injectivity of $\chi_q$ is a fusion product, that is to say the constant structures on the basis of simple representations are positive integers.
\end{thm}
A natural question is to "categorify" $*$, to construct a tensor category of representations of $\U_q(\hat{\Glie})$ corresponding to $*$. In \cite{h5} such a category is proposed by using the Drinfeld coproduct (in the particular type $A_1^{(1)}$, another procedure is discussed in \cite{m3}). In the rest of this section we restrict to the cases where the Drinfeld coproduct is known to be compatible with Drinfeld-Serre relations, that is to say when $(i\neq j\Rightarrow C_{i,j}C_{j,i}\leq 3)$ (this includes almost all quantum toroidal algebras; when the condition is not satisfied analog results hold for the algebra without Drinfeld-Serre relations \cite{h5}). The Drinfeld coproduct is defined in a completion, and so can not be used in the usual way to construct tensor products of representations (we get divergence problems). To overcome this problem, we use a deformation/renormalization process. The strategy is analog to the construction of crystal bases by Kashiwara. The deformed Drinfeld coproduct is given by the following formulas on generators of $\U_q(\hat{\Glie}) ((u))$:
$$\Delta_u(x_{i,r}^+)= x_{i,r}^+\otimes 1 + \underset{l\geq 0}{\sum} u^{r+l} (\phi_{i,-l}^-\otimes x_{i,r+l}^+),$$
$$\Delta_u(x_{i,r}^-)=u^r (1\otimes x_{i,r}^-) + \underset{l\geq 0}{\sum} u^l (x_{i,r-l}^-\otimes \phi_{i,l}^+),$$
$$\Delta_u(\phi_{i,\pm m}^{\pm})= \underset{0\leq l\leq m}{\sum}u^{\pm l} (\phi_{i,\pm (m-l)}^{\pm}\otimes\phi_{i,\pm l}^{\pm})\text{ , }\Delta_u(k_h)=k_h\otimes k_h,$$
where $u$ is a formal deformation parameter.

\begin{rem} Let us make precise the claim that the non deformed Drinfeld coproduct makes sense in a completion. $\U_q(\hat{\Glie})$ has a grading defined by $\text{deg}(x_{i,m}^\pm) = m$, $\text{deg}(\phi_{i,\pm m}^\pm) = \pm m$, $\text{deg}(k_h) = 0$. Let $\U_q(\hat{\Glie})\hat{\otimes}\U_q(\hat{\Glie})$ be the topological completion of $\U_q(\hat{\Glie})\otimes \U_q(\hat{\Glie})$ with respect to the degree of the second term, that is to say the $\U_q(\hat{\Glie})$-module of sums of the form $\sum_{R\geq 0} a_R\otimes b_R$ where $a_R,b_R\in\U_q(\hat{\Glie})$, $b_R$ is homogoneous and $\text{deg}(b_R)$ has limit $+\infty$ when $R\rightarrow +\infty$. Then $\U_q(\hat{\Glie})\hat{\otimes}\U_q(\hat{\Glie})$ has a natural algebra structure and the specializations at $u = 1$ of the formulas above for $\Delta_u(x_{i,r}^+)$, $\Delta_u(x_{i,r}^-)$, $\Delta_u(\phi_{i,\pm m}^{\pm})$, $\Delta_u(k_h)$ make sense in $\U_q(\hat{\Glie})\hat{\otimes}\U_q(\hat{\Glie})$. We get a well-defined algebra morphism $\Delta : \U_q(\hat{\Glie}) \rightarrow \U_q(\hat{\Glie})\hat{\otimes}\U_q(\hat{\Glie})$. Let us check the coassociativity. Let $\U_q(\hat{\Glie})^{\hat{\otimes}3}$ be the $\U_q(\hat{\Glie})$-module of sums of the form $\sum_{R\geq 0} a_R \otimes b_R\otimes c_R$  where $b_R,c_R$ are homogeneous and $2\text{deg}(c_R) + \text{deg}(b_R)$ has limit $+\infty$ when $R\rightarrow + \infty$. It has an algebra structure. From the formulas in \cite[Section 3.1]{h5} we see that $(\text{Id}\otimes \Delta)\circ \Delta$ and $(\Delta\otimes \text{Id})\circ\Delta$ make sense as algebra morphisms $\U_q(\hat{\Glie})\rightarrow \U_q(\hat{\Glie})^{\hat{\otimes}3}$. Moreover from the relation $(\text{Id}\otimes \Delta_u)\circ \Delta_u = (\Delta_u\otimes \text{Id})\circ \Delta_{u^2}$ in \cite[Lemma 3.4]{h5}, these two algebra morphisms are equal.\end{rem}

In \cite{h5} a category $\text{Mod}$ of $\U_q(\hat{\Glie})\otimes \CC(u)$-modules satisfying a certain quasipolynomiality property is considered, and it is proved that $\Delta_u$ defines a tensor category structure on $\text{Mod}$. The associativity axiom is proved from the following twisted coassociativity of $\Delta_u$ ($r,r'\in\ZZ$):
$$(\text{Id}\otimes\Delta_{u^{r'}})\circ\Delta_{u^r}=(\Delta_{u^r}\otimes\text{Id})\circ\Delta_{u^{r+r'}}.$$
We do not give in this review the detail of the construction for the whole category $\text{Mod}$ (we refer to \cite{h5}), but we explain what is the most relevant for the considered applications : let $V$ and $V'$ be two $l$-highest weight integrable representations of respective $l$-highest weight $v$ and $v'$. $\Delta_u$ defines a structure of $\U_q(\hat{\Glie})\otimes\CC((u))$-module on $(V\otimes V')\otimes \CC((u))$. In fact $(V\otimes V')\otimes \CC(u)$ is stable for the action of $\U_q(\hat{\Glie})\otimes\CC(u)$, and the fusion module is by definition :
$$V *_f V' = (\U_q(\hat{\Glie})\otimes\CC(u)).(v\otimes v')/(u-1)(\U_q(\hat{\Glie})\otimes\CC(u))(v\otimes v').$$
\begin{thm}\cite{h2} $V *_f V'$ is an $l$-highest weight $\U_q(\hat{\Glie})$-module and 
$$[V*_f V'] = [V]*[V'],$$
where $[.]$ denotes the image in the Grothendieck group.
\end{thm}
Remark : the last result can be interpreted as a cyclicity property as this gives a module $V*_f V'$ of $l$-highest weight. In particular for quantum affine algebras, $V*_f V'$ is not the usual tensor product. Indeed for the usual tensor product of quantum affine algebras, certain tensor products of $l$-highest weight modules are $l$-highest weight (see \cite[Theorem 4]{c} and \cite[Theorem 9.1]{kas0}) but not all. For example for $\mathcal{U}_q(\hat{sl_2})$, $V(a)\otimes V(aq^2)$ is not of $l$-highest weight, but $V(a) * V(aq^2)\simeq V(aq^2)* V(a)\simeq V(aq^2)\otimes V(a)$ are of $l$-highest weight (for the other choice of the coproduct, $V(a)\otimes V(aq^2)$ is of $l$-highest weight). 

By analogy to the case of quantum affine algebras, we expect the existence of maximal $l$-highest weight integrable representations (Weyl modules). By the above construction, such modules should admit a fusion product of fundamental representations as a subquotient. We conjecture that the Weyl modules exist and are isomorphic to a fusion product of fundamental representations.

\subsection{Examples : $T$-system for Kirillov-Reshetikhin modules}

The Kirillov-Reshetikhin modules satisfy the $T$-system :

\begin{thm}\cite{h5}\label{th5} We have an exact sequence :
$$0\rightarrow S_{k,a}^{(i)}\rightarrow W_{k,a}^{(i)}*_f W_{k,aq_i^2}^{(i)}\rightarrow W_{k+1,a}^{(i)} *_f W_{k-1,aq_i^2}^{(i)}\rightarrow 0,$$
where :
$$S_{k,a}^{(i)}=L_{\nu_{k,a}^{(i)}}*_f({*_f}_{(j,l)\in A_i}W_{K(j,l),aq_j^{-(2l-1)/C_{i,j}}}^{(j)}),$$
and :
$$\nu_{k,a}^{(i)}=k(\Lambda_i-\alpha_i)-{\sum}_{(j,l)\in A_i}K(j,l)\Lambda_j,$$
$$A_i = \{(j,l)\in I\times\ZZ|C_{j,i}<0\text{ , }1\leq l\leq -C_{i,j}\},$$
$$K(j,l) = -C_{j,i}+E(r_i(k-l)/r_j),$$
and $E(m)\in\ZZ$ denotes the integral part of $m\in\RR$.
\end{thm}
\noindent For quantum affine algebras, an analog result holds with the usual tensor product, and was proved in \cite{Nab, Nad} for simply-laced cases and in \cite{hc} for non simply-laced cases. It is the crucial point for the proof of the Kirillov-Reshetikhin conjecture.

\begin{rem} The exact sequence of Theorem \ref{th5} is not splitted as all involved modules are $l$-highest weight, and so $W_{k,a}^{(i)}*_f W_{k,aq_i^2}^{(i)}$ is not decomposable.\end{rem}

\begin{rem}\label{vgr} Viewed in the Grothendieck group, the exact sequence provides a combinatorial inductive system called $T$-system. Viewed in $\mathcal{E}$, we get a solution which is polynomial in the monomials $m_\gamma$ with positive coefficients : this is analog to the Fomin-Zelevinsky Laurent Phenomena positivity properties for analog systems ($Y$-systems) \cite{fz1, fz2} as noticed in \cite{h8}.\end{rem}

\noindent In the case of the quantum affine algebra $\U_q(\hat{sl}_{n+1})$-modules ($n\geq 2$) we have :
$$0\rightarrow W_{1,aq}^{(2)} \rightarrow W_{1,a}^{(1)}\otimes W_{1,aq^2}^{(1)}\rightarrow W_{2,a}^{(1)}\rightarrow 0.$$
For $\epsilon = q^2$, we get via Schur-Weyl duality the exact sequence of $\tilde{\mathcal{H}}_2(\epsilon)$-modules :
$$0\rightarrow V_{(a,a\epsilon)}\rightarrow M_{(a,a\epsilon)}\rightarrow V_{(a),(a\epsilon)}\rightarrow 0.$$
Here $M_{(a,a\epsilon)}$ is not semi-simple and $\CC(\sigma - \epsilon)$ is a submodule isomorphic to $V_{(a,a\epsilon)}$. Let us study $M_{(a\epsilon,a)}$. It has a submodule $\CC.(1+\sigma)$ isomorphic to $V_{(a),(a\epsilon)}$ and it is not semi-simple (the eigenspaces of $z_1$ are not submodules). We have :
$$0\rightarrow V_{(a),(a\epsilon)}\rightarrow M_{(a\epsilon,a)}\rightarrow V_{(a,a\epsilon)}\rightarrow 0.$$
It follows from this discussion that $W_{1,aq^2}^{(1)}\otimes W_{1,a}^{(1)}$ is not semi-simple and that we have a non splitted exact sequence : 
$$0\rightarrow W_{2,a}^{(1)} \rightarrow W_{1,aq^2}^{(1)}\otimes W_{1,a}^{(1)}\rightarrow W_{1,aq}^{(2)}\rightarrow 0.$$
%This example confirms that the condition $n\geq l$ is crucial to have the equivalence of category in Theorem \ref{cp2t}.

\subsection{Additional results}

A particular class of simple representations (called special modules, notion defined in \cite{Nab}) is of particular importance as a generalization of an algorithm of Frenkel-Mukhin \cite{Fre2} gives the $q$-character of these representations (see the discussion at the end of \cite{npc}) :

\begin{defi} An integrable representation $V$ in the category $\mathcal{O}$ is said to special if it has a unique $l$-weight satisfying condition of Theorem \ref{cond} and if its $l$-weight space is of dimension $1$.
\end{defi}
The Kirillov-Reshetikhin modules are special for quantum affine algebras \cite{Nab, Nad, hc} and in general \cite{h5}. But there are much more special representations, as for example a large class of minimal affinizations of quantum affine algebras \cite{h4} (generalizations of Kirillov-Reshetikhin modules interesting from the physical point of view, as stressed in \cite[Remark 4.2]{frere} and in the introduction of \cite{Chari2}). As an illustration let us look at an example for more general quantum affinizations. Suppose that $\Glie$ is of type $B_{n,p}$ ($n\geq 2$, $p\geq 2$) of Cartan matrix $C_{i,j} = 2\delta_{i,j} - \delta_{j,i+1} - \delta_{j,i-1} -(p-1) \delta_{i,n}\delta_{j,n-1}$ (in fact $B_{n,1} = A_n$, $B_{n,2} = B_n$, $B_{2,3} = G_2$). Let $\lambda\in P^+$. For $i\in I$ let $\lambda_i=\lambda(\alpha_i^{\vee})$ and for $i < n$ consider analogs of the Chari-Pressley coefficients classifying the minimal affinizations :
$$c_i(\lambda)=q^{r_i\lambda_i+r_{i+1}\lambda_{i+1}+r_{i}-C_{i,i+1}-1}.$$ 
Let us denote by $\gamma_{k,a}^{(i)}$ the $l$-highest weight of $W_{k,a}^{(i)}$. 
\begin{thm}\cite{h4} For $\gamma$ of the form $\gamma = \prod_{i\in I} \gamma_{\lambda_i,a_i}^{(i)}$ with $(a_i)_{i\in I}\in(\CC^*)^I$, $L(\gamma)$ is special if :
$$\forall i<j\in I\text{ , }a_j/a_i=\prod_{i\leq s <j}c_s'(\lambda).$$
\end{thm}
The small property, stronger than the special property, is related to the geometric smallness in the sense of Borho-MacPherson. To define this property, we need the combinatorial description of $q$-characters. By generalizing the result in \cite{Fre, Naams}, it is proved \cite{h6} that for general quantum affinization an $l$-weight $\gamma$ satisfies the property (the main point is that $\hat{\mathcal{U}}_i$ is isomorphic to a quantum affine algebra of type $sl_2$):
$$\gamma(\phi_i^\pm (z)) = q_i^{\text{deg}(Q_i)-\text{deg}(R_i)}\frac{Q_i(zq_i^{-1})R_i(zq_i)}{Q_i(zq_i)R_i(zq_i^{-1})},$$
where $Q_i(z),R_i(z)\in\CC[z]$ have constant term equal to $1$. So we can define the monomial $m_\gamma = \prod_{i\in I}Y_{i,a}^{\mu_{i,a} - \nu_{i,a}}$ where
$$Q_i(z)={\prod}_{a\in\mathbb{C}^*}(1-za)^{\mu_{i,a}}\text{ and }R_i(z)={\prod}_{a\in\mathbb{C}^*}(1-za)^{\nu_{i,a}}.$$ 
The variables $Y_{i,a}^{\pm 1}$ originally appeared from the $R$-matrix as element of the completed Cartan subalgebra (see the original paper \cite{Fre}, and \cite{h2} for general quantum affinizations). $e(\gamma)$ is also denoted by $m_\gamma$. For $i\in I, a\in\CC^*$, consider the monomial :
$$A_{i,a}=Y_{i,aq_i^{-1}}Y_{i,aq_i}{\prod}_{\{(j,k)|C_{j,i}\leq -1\text{ and }k\in\{C_{j,i}+1,C_{j,i}+3,...,-C_{j,i}-1\}\}}Y_{j,aq^k}^{-1}.$$
We write $\gamma\leq \gamma'$ if $m_\gamma$ is the product of $m_{\gamma'}$ by some $A_{i,a}^{-1}$. We have the following generalization \cite{h3} of a result of \cite{Fre2, Naams} :
\begin{thm} The $l$-weights $\gamma'$ of an integrable $L(\gamma)$ satisfy $\gamma\leq \gamma'$.
\end{thm}
Borho and MacPherson introduced \cite[Section 1.1]{bm} remarkable geometric properties (smallness and semi-smallness) for a proper algebraic map $\pi : Z \rightarrow X$ where $Z, X$ are irreducible complex algebraic varieties. This geometric situation is of particular interest as the Beilinson-Bernstein-Deligne-Gabber decomposition Theorem is simplified \cite[Section 1.5]{bm}. The Springer resolution is a fundamental example of a semi-small morphism \cite{bm}. Nakajima defined the intensively studied quiver varieties. They depend on a quiver $Q$ and come with a resolution which is semi-small for a finite Dynkin diagram (see \cite[Section 5.2]{npc}). The graded version of quiver varieties are also of particular importance, for example for their relations with quantum affine algebras \cite{Nab}. They also come with resolutions. A natural problem addressed in \cite{h6} is to study the small property of these resolutions (\cite[Conjecture 10.4]{Nab}). In fact the geometric small property can be translated in terms of small module of quantum affine algebras. We suppose that $\Glie$ is simply-laced. 
\begin{defi} $L(\gamma)$ is small if and only if for all $\gamma'\leq \gamma$, $L(\gamma')$ is special.\end{defi}
\noindent The above definition is a purely representation theoretical characterization of the property which was noticed in \cite{h6} by refining a proof of \cite{Nab}.

\begin{defi} A node $i\in \{1,\cdots ,n\}$ is said to be extremal (resp. special) if there is a unique $j\in I$ (resp. three distinct $j,k,l\in I$) such that $C_{i,j} < 0$ (resp. $C_{i,j}<0$, $C_{i,k}<0$ and $C_{i,l}<0$). We denote by $d_i$ the minimal $d\in\NN\cup\{+\infty\}$ such that there are distinct $i_1,\cdots,i_d\in I$ satisfying $C_{i_j,i_{j+1}} < 0$ and $i_d$ is special.\end{defi}

We get for general simply-laced quantum affinizations :

\begin{thm}\cite{h6} For $k > 0, i\in I, a\in\CC^*$, $W_{k,a}^{(i)}$ is small if and only if 
$$k\leq 2\text{ or }(i\text{ is extremal and }k\leq d_i+1).$$\end{thm} 

\section{New results and further developments}\label{trois}

In this section we announce new results and explain possible further developments illustrated with examples.

\subsection{Relation to $\U_q(\hat{sl}_\infty)$ and character formulas}

We consider the quantum toroidal algebra $\U_q(sl_{n+1}^{tor})$ ($n\geq 2$). We identify $I$ with $\ZZ/(n+1)\ZZ$. We have exact sequences of Theorem \ref{th5} : 
$$0\rightarrow W_{k,aq}^{(i + 1)} *_f W_{k,aq}^{(i-1)}\rightarrow W_{k,a}^{(i)}*_f W_{k,aq^2}^{(i)}\rightarrow W_{k+1,a}^{(i)} *_f W_{k-1,aq_i^2}^{(i)}\rightarrow 0.$$
So the $\mathcal{W}_{k,a}^{(i)} = \chi_q(W_{k,a}^{(i)})$ satisfy the corresponding $T$-system :
$$\mathcal{W}_{k,a}^{(i)})\mathcal{W}_{k,aq^2}^{(i)} = \mathcal{W}_{k,aq}^{(i + 1)}\mathcal{W}_{k,aq}^{(i-1)} + \mathcal{W}_{k+1,a}^{(i)}\mathcal{W}_{k-1,aq^2}^{(i)}.$$
Consider the ring morphism 
$$R:\ZZ[Y_{i,a}^{\pm 1}]_{i\in I,a\in\CC^*}\rightarrow \ZZ[Y_{i,a}^{\pm 1}]_{i\in I, a\in\CC^*}$$ 
defined by $R(Y_{i,a}) = Y_{i+1,a}$. By symmetry we have $\mathcal{W}_{k,a}^{(i)} = R^i(\mathcal{W}_{k,a}^{(0)})$ for any $k\geq 0$, $a\in\CC^*$.

\noindent For $a\in\CC^*$ and $i\in \ZZ$ let $\ffbox{i}_a = Y_{[i-1],aq^i}^{-1} Y_{[i],aq^{i-1}}$ where $[j]\in I$ is the class of $j\in\ZZ$. Let $\mathcal{T}$ be the set of tableaux $(T_{i,j})_{i\leq 0, 1\leq j\leq k}$ such that 

$\forall i\leq 0, 1\leq j\leq k$, $T_{i,j}\in \ZZ$,

$\forall i\leq 0, 1\leq j\leq k-1$, $T_{i,j}\leq T_{i,j+1}$,

$\forall i\leq 0, 1\leq j\leq k$, $T_{i,j} < T_{i+1,j}$,

$\forall 1\leq j\leq k$, $T_{i,j} = i$ for $i << 0$.

\begin{thm}\label{nexp} We have the following explicit formula :
$$\mathcal{W}_{k,aq}^{(\ell)} = R^\ell(\sum_{T\in \mathcal{T}} \prod_{i\leq 0,1\leq j\leq k} \ffbox{T_{i,j}}_{aq^{2(j-i)}}).$$
\end{thm}

\noindent Remark : in the case of fundamental representations ($k = 1$) explicit formulas have been proved in \cite{nag} as a special case of $q$-Fock space representations.

\noindent The details of the direct proof of the Theorem will be written in a separate publication devoted to the general study of the quantum affinization Kac-Moody algebra $\U_q(\hat{sl}_\infty)$ with infinite Cartan matrix $(C_{i,j})_{i,j\in\ZZ}$ where 
$$C_{i,j} = 2\delta_{i,j} - \delta_{i,j+1} - \delta_{i,j-1}$$ 
(the algebra is defined as in Definition \ref{defiqaf} but $i\in \ZZ$). The corresponding infinite Dynkin diagram admits a non trivial automorphism $\tau_m : i\mapsto i+m$ for $m\geq 2$ which gives a twisted Dynkin diagram of type $A_{m-1}^{(1)}$. So we can use technics developed in \cite{h8} to study twisted quantum affine algebras to get the formula of Theorem \ref{nexp}. 

\noindent Although for quantum affine algebras of type $A,B,C$ and for the double deformed quantum toroidal algebras of type $A$ (of Remark \ref{tp}) \cite{nag} the $l$-weight spaces of fundamental representations are of dimension $1$, this is not true for $\U_q(sl_{n+1}^{tor})$. For example we have from the Frenkel-Mukhin algorithm the upper part of the $q$-character, presented bellow as a graph (see \cite{Fre} for details on the notion of graph associated with a $q$-character). 

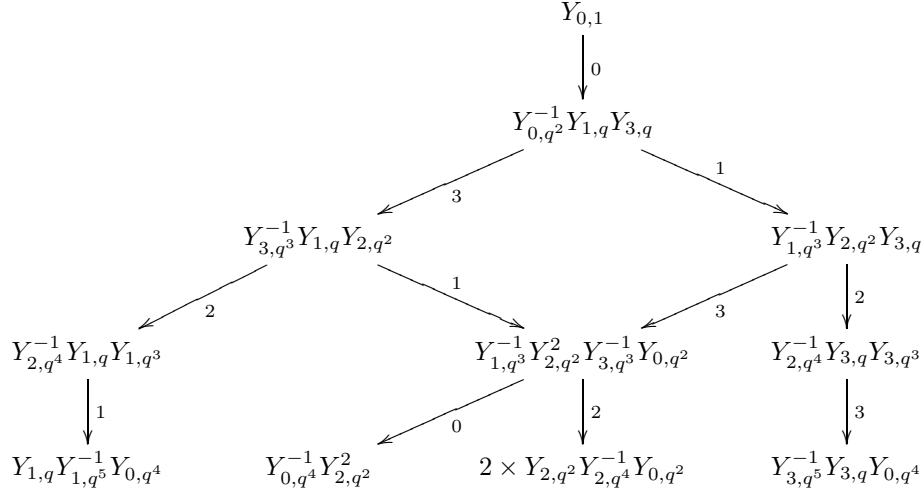
\begin{figure}[htbp]\label{fid}
\begin{center}
\begin{equation*}
\xymatrix{&&Y_{0,1} \ar[d]^0&
\\&&Y_{0,q^2}^{-1}Y_{1,q}Y_{3,q}\ar[ld]^3\ar[rd]^1&
\\&Y_{3,q^3}^{-1}Y_{1,q}Y_{2,q^2}\ar[rd]^1\ar[dl]^2&&Y_{1,q^3}^{-1}Y_{2,q^2}Y_{3,q}\ar[ld]^3\ar[d]^2
\\Y_{2,q^4}^{-1}Y_{1,q}Y_{1,q^3}\ar[d]^1&&Y_{1,q^3}^{-1}Y_{2,q^2}^2Y_{3,q^3}^{-1}Y_{0,q^2}\ar[d]^2\ar[ld]^0&\ar[d]^3Y_{2,q^4}^{-1}Y_{3,q}Y_{3,q^3}
\\Y_{1,q}Y_{1,q^5}^{-1}Y_{0,q^4}&Y_{0,q^4}^{-1}Y_{2,q^2}^2&2\times Y_{2,q^2}Y_{2,q^4}^{-1}Y_{0,q^2}&Y_{3,q^5}^{-1}Y_{3,q}Y_{0,q^4}}
\end{equation*}
\caption{Type $A_3^{(1)}$ : a part of the graph of $\chi_q(V_0(1))$}
\end{center}
\end{figure}

\begin{rem} Note that as $\text{dim}((L(\Lambda_0))_{\Lambda_0 - \alpha_0-\alpha_3-\alpha_1-\alpha_2}) = 1$, $V_0(1)$ is not simple as a $\U_q^h(\hat{sl}_4)$-module (this is not a situation analog to evaluation representations of quantum affine algebras).\end{rem}

Let us give some additional comments on representations of $\U_q(\hat{sl}_\infty)$. We can define Kirillov-Reshetikhin modules $\tilde{W}_{k,a}^{(i)}$ ($k\geq 0,a\in\CC^*,i\in\ZZ$) of $\U_q(\hat{sl}_\infty)$ by using Drinfeld polynomials as for $\U_q(\hat{sl}_{n+1})$. These modules are not finite dimensional, but by using the same proof as for quantum affinizations \cite{h2}, they are integrable. Then, as for quantum affinizations we get :
\begin{thm} The $\tilde{W}_{k,a}^{(i)}$ satisfy the $T$-system for any $k\geq 0,a\in\CC^*,i\in\ZZ$ :
$$\chi_q(\tilde{W}_{k,a}^{(i)})\chi_q(\tilde{W}_{k,aq^2}^{(i)}) = \chi_q(\tilde{W}_{k,aq}^{(i + 1)})\chi_q(\tilde{W}_{k,aq}^{(i-1)}) + \chi_q(\tilde{W}_{k+1,a}^{(i)})\chi_q(\tilde{W}_{k-1,aq^2}^{(i)}).$$
\end{thm}
In fact by putting $T_{i,k,t} = \chi_q(\tilde{W}_{k,aq^{t+1-k}}^{(i)})$ we get the octahedron recurrence (see \cite{rr} and more recently \cite{fz1, ktw, hk} for related works) :
$$T_{i,k,t-1}T_{i,k,t+1} = T_{i+1, k, t}     T_{i-1,k,t} + T_{i,k+1,t}T_{i,k-1,t}.$$
In particular the comments of Remark \ref{vgr} hold for this relation.

\subsection{Extremal loop modules and finite dimensional representations at roots of unity}

In this section we introduce new class of representations of quantum toroidal algebras, and we give motivations and examples to study them. First let us recall the definition of extremal weight modules. The Weyl group $W$ is the subgroup of $GL(\Hlie)$ generated by the simple reflection $s_i$ defined by $s_i(\lambda)=\lambda-\lambda(\alpha_i^{\vee})\alpha_i$ ($i\in I$).

\begin{defi}\cite{kasm1, kas0} For $V$ an integrable $\U_q(\Glie)$-module and $\lambda\in P$, a vector $v\in V_{\lambda}$ is called  extremal of weight $\lambda$ if there are vectors $\{v_w\}_{w\in W}$ such that $v_{\text{Id}}=v$ and for any $i\in I$ :
$$x_i^{\pm} v_w=0 \text{ if } \pm w(\lambda)(\alpha_i^{\vee})\geq 0
\text{ and }(x_i^{\mp})^{\pm (w(\lambda)(\alpha_i^{\vee}))} v_w=v_{s_i(w)}.$$
For $\lambda\in P$, the {\it extremal weight module $V(\lambda)$ of extremal weight $\lambda$\/} is the $\U_q(\Glie)$-module generated by a
vector $v_{\lambda}$ with the defining relations that $v_{\lambda}$ is extremal of weight $\lambda$.\end{defi}
If $\lambda\in P^+$, $V(\lambda)$ is the
simple module of highest weight $\lambda$, but in general $V(\lambda)$ is not in $\mathcal{O}$.
\begin{thm}\cite{kasm1}
For $\lambda\in P$, the module $V(\lambda)$ is
integrable and has a crystal basis $\mathcal{B}(\lambda)$.
\end{thm}
 By a monomial we mean a monomial of $\ZZ[Y_{i,q^l}^{\pm}]_{1\leq i\leq n, l\in\ZZ}$. In \cite{hn} a monomial realization of extremal crystals is proved, by extending the results in \cite{Nac, kas} obtained for crystals of finite type. The monomial crystal of a crystal of finite type is related to the $q$-character of a finite dimensional representation of a quantum affine algebra \cite{Nac}. We propose a relation between level $0$ monomial crystal \cite{hn} and some integrable representations of quantum toroidal algebras (which are not in $\mathcal{O}$). Let us look at an example in the most simple situation $A_3^{(1)}$ ($A_1^{(1)}$ is not simply-laced, and $A_2^{(1)}$ has an odd cycle and so is not included in the study of \cite{hn}). The crystal $\mathcal{B}(\Lambda_1 - \Lambda_0)$ can be realized as :
$$\cdots\overset{0}{\rightarrow} Y_{1,q^0}Y_{0,q}^{-1}\overset{1}{\rightarrow} Y_{2,q}Y_{1,q^2}^{-1}\overset{2}{\rightarrow} Y_{3,q^2}Y_{2,q^3}^{-1}
\overset{3}{\rightarrow}Y_{0,q^3}Y_{3,q^4}^{-1} \overset{0}{\rightarrow}Y_{1,q^4}Y_{0,q^5}^{-1}\overset{1}{\rightarrow}\cdots$$
\begin{prop}\label{exa} There is an integrable representation $V$ of $\U_q(sl_4^{tor})$ such that :
$$\chi_q(V) = \sum_{p\in\ZZ} (Y_{1,q^{4p}}Y_{0,q^{4p+1}}^{-1} + Y_{2,q^{4p+1}}Y_{1,q^{4p+2}}^{-1} + Y_{3,q^{4p+2}}Y_{2,q^{4p+3}}^{-1} + Y_{0,q^{4p+3}}Y_{3,q^{4p+4}}^{-1}).$$
\end{prop}

\demo
Consider the vector space $V = \bigoplus_{p\in\ZZ} V_p$ where 
$$V_p = \oplus_{a\in\{1,2,3,4\}}\CC v_{a,p}.$$ 
Let us define the action of $\U_q(sl_4^{tor})$ on $V$ by ($r\in\ZZ$, $m > 0$, $a,b\in\{0,\cdots,3\}$, $p\in\ZZ$) :
$$x_{a,r}^+.v_{a+b,p} = \delta_{b,1}q^{r(4p+a-1)}v_{a,p},$$
$$x_{a,r}^-.v_{a+b,p} = \delta_{b,0} q^{r(4p + a-1)}v_{a+1,p},$$
$$\phi_{a,\pm m}^{\pm}.v_{a+b,p} = \pm (\delta_{b,0} - \delta_{b,1}) (q - q^{-1})q^{\pm m(4p + a-1) }v_{a+b,p},$$
$$k_a^{\pm}.v_{a+b,p} = q^{\pm(\delta_{b,0} - \delta_{b,1})}v_{a+b,p},$$  
where we denote $v_{0,p} = v_{4,p-1}$, $v_{5,p} = v_{1,p+1}$, $v_{6,p} = v_{2,p+1}$. Let us check that it defines a representation of $\U_q(sl_4^{tor})$.  For each $a\in\{0,1,2,3\}$, we have a well defined action of $\hat{\U}_a$. So it suffices to check the relation involving adjacent nodes, $(0,1)$, $(1,2)$, $(2,3)$, $(3,4)$. For example for $(2,3)$, we can see that $V_p$ is stable for $\hat{\U}_{1,2}$ and is isomorphic to the direct sum of $(\CC v_{1,p})$ (trivial representation) and $\CC v_{2,p} \oplus \CC v_{3,p} \oplus \CC v_{4,p}$ (a fundamental representation for the node $2$).\qed

\noindent This representation can be considered as an analog $V(Y_{1,1}Y_{0,q}^{-1})$ of extremal module for quantum toroidal algebra : $Y_{1,1}Y_{0,q}^{-1}$ is the $l$-weight of the extremal vector $v_{1,0}$. For $0\leq i\leq n$, let $\U_q^{v,i}(\hat{\Glie})$ be the subalgebra of $\U_q(\Glie^{tor})$ generated by the $x_{j,m}^{\pm}, k_j^{\pm 1},h_{j,r}$ for $j\neq i$. $\U_q^{v,i}(\hat{\Glie})$ is isomorphic to a quantum affine algebra and is equal to the vertical subalgebra when $i$ is the additional node of the Dynkin diagram. In general we propose the following definition : 

\begin{defi} An extremal loop module of $\U_q(\Glie^{tor})$ is an integrable representation $V$ such that there is $v\in V$ $l$-weight vector satisfying 

(1) $\U_q(\Glie^{tor}).v = V$, 

(2) $v$ is extremal for $\U_q^h(\hat{\Glie})$, 

(3) $\forall w\in V$, $\forall 0\leq i\leq n$, $\U_q^{v,i}(\hat{\Glie}).w$ is finite dimensional.
\end{defi}
In another paper such integrable representations will be studied (as well as the analogy with \cite{cg, v} for quantum affine algebras). 

Let us go back to our example. As a $\U_q^v(\hat{sl}_4)$-module $V$ is not simple and is semi-simple equal to a direct sum of simple finite dimensional representations. As a $\U_q^h(\hat{sl}_4)$-module it is extremal with a non trivial automorphism satisfying $v_{a,p}\mapsto v_{a,p+1}$. But the $\U_q(sl_4^{tor})$-structure of $V$ has very different properties than its $\U_q^v(\hat{sl}_4)$ and $\U_q^h(\hat{sl}_4)$-structures as it is simple without non trivial automorphism (and so it is {\it not} analog to evaluation representations of quantum affine algebras).

\begin{rem} We can expect that for a large class of such representations the level is finite, and the representation can be interpreted via Schur-Weyl duality. In our example the character of $V_p$ as a $\U_q(sl_4)$-module is $y_1 + y_2y_1^{-1} + y_3y_2^{-1} + y_3^{-1}$ and so $V$ is of level $1$. We get a representation $M$ of $\mathbb{H}_1(q^2,q^4)$ (the twisted polynomials ring with two variable $X_1,Y_1$ and relation $X_1Y_1 = q^4Y_1 X_1$). We have $M = \bigoplus_{p\in\ZZ} \CC m_p$ where $X_1m_p = m_{p-1}$ and $Y_1m_p = q^{4p}m_{p+1}$.\end{rem}

Kashiwara \cite{kas0} proved that certain fundamental finite dimensional representations of quantum affine algebras can be realized as a quotient of a level 0 fundamental extremal representation, and so have a crystal basis. The quotient comes from an automorphism of fundamental level $0$ extremal modules. In many case \cite{hn}, this periodicity property can be read on the monomial model as a shift on the index by some $q^l$. This is related to the existence of finite dimensional representations of quantum toroidal algebras at roots of unity. Moreover we expect to get finite level representations and so via Schur-Weyl duality to get finite dimensional representations of DAHA at roots of unity (note that in the case of the $q$-Fock space, the Schur-Weyl duality is used in the reversed way).

\noindent In the above studied example, for $q=i$ a primitive $4$-root of unity, we get a periodic sequence of monomials :
$$\cdots\overset{0}{\rightarrow} Y_{1,1}Y_{0,i}^{-1}\overset{1}{\rightarrow} Y_{2,i}Y_{1,-1}^{-1}\overset{2}{\rightarrow} Y_{3,-1}Y_{2,-i}^{-1}
\overset{3}{\rightarrow}Y_{0,-i}Y_{3,1}^{-1} \overset{0}{\rightarrow}Y_{1,1}Y_{0,i}^{-1}\overset{1}{\rightarrow}\cdots
$$
Let $L\geq 1$ and $\epsilon$ a $4L$-primitive root of unity. Then we also get a periodic sequence when $q$ is equal to $\epsilon$. Let us consider $\U_\epsilon(sl_4^{tor})$ the algebra defined as $\U_q(sl_4^{tor})$ with $\epsilon$ instead of $q$ (we do not include divided powers in the algebra).

\begin{prop} There is a $4L$-dimensional representation $\tilde{V}$ of $\U_\epsilon(sl_4^{tor})$ such that :
$$\chi_\epsilon(\tilde{V}) = \sum_{1\leq p\leq L} (Y_{1,\epsilon^{4p}}Y_{0,\epsilon^{4p+1}}^{-1} + Y_{2,\epsilon^{4p+1}}Y_{1,\epsilon^{4p+2}}^{-1} + Y_{3,\epsilon^{4p+2}}Y_{2,\epsilon^{4p+3}}^{-1} + Y_{0,\epsilon^{4p+3}}Y_{3,\epsilon^{4p+4}}^{-1}).$$
\end{prop}

\demo It suffices to specialize at $q = \epsilon$ the representation of Proposition \ref{exa} and to identify $v_{a,p}$ with $v_{a,p+L}$ for any $p$.\qed

\begin{rem} In this example the corresponding representation of $\mathbb{H}_1(\epsilon^2,\epsilon^4)$ (the algebra with generators $X,Y$ and relation $XY = \epsilon^4 YX$) is $W$ of dimension $L$ with basis $\{w_1,\cdots,w_L\}$ where $Xw_i=\epsilon^{4i}w_i$ and $Yw_i = w_{i+1}$ (we denote $w_{L+1} = w_1$). Note that for $\{m_1,\cdots,m_L\}$ the basis of $W$ such that $w_i = \sum_{1\leq j\leq L}\epsilon^{4ij}m_j$, we have $X m_j = m_{j-1}$ and $Y m_j = \epsilon^{4j} m_j$. The eigenvalues of $X$ and of $Y$ are $\{1,\epsilon^4,\cdots,\epsilon^{4(L-1)}\}$.\end{rem}

\begin{rem} In the case $L = 1$ we have a $4$ dimensional representation of $\U_i(sl_4^{tor})$ with $i$-character
$$\chi_i(\tilde{V}) = Y_{1,1}Y_{0,i}^{-1} + Y_{2,i}Y_{1,-1}^{-1} + Y_{3,-1}Y_{2,-i}^{-1} + Y_{0,-i}Y_{3,1}^{-1}.$$ 
More explicitly $\tilde{V} = \oplus_{a\in\ZZ/4\ZZ}\CC v_a$ and the action is given by ($r\in\ZZ$, $m > 0$, $a\in\{0,\cdots,3\}$, $a'\in\ZZ/4\ZZ$) :
$$x_{a,r}^+.v_{a'} = \delta_{a',[a+1]}i^{r(a-1)}v_{a' - 1},$$
$$x_{a,r}^-.v_{a'} = \delta_{a',[a]} i^{r(a-1)}v_{a'+1},$$
$$\phi_{a,\pm m}^{\pm}.v_{a'} = \pm 2(\delta_{a',[a]} - \delta_{a',[a+1]}) i^{1 \pm m(a-1) }v_{a'},$$
$$k_a^{\pm}.v_{a'} = i^{\pm(\delta_{a',[a]} - \delta_{a',[a+1]})}v_{a'}.$$  
The structures of $\U_i^h(\hat{sl_4})$ and $\U_i^v(\hat{sl_4})$-module of $\tilde{V}$ are both isomorphic to a simple $4$-dimensional fundamental representation. In this simple example the corresponding representation of $\mathbb{H}_1(-1,1)\simeq \CC[X_1,Y_1]$ is of dimension $1$ where $X_1,Y_1$ act as the identity.\end{rem}

\noindent These representations for general type will be studied systematically in another paper.

\subsection{Other possible developments}

In this last section we explain other possibly promising directions for further research. A first general observation is that in the view of the huge number of very interesting results which are known or under development for Cherednik algebras (as for example in this volume), we can expect to have analog results, developments and applications for quantum toroidal algebras. Let us give some examples and more precise directions :

- Many combinatorial results and explicit formulas are known for $q$-characters of Kirillov-Reshetikhin modules, in particular for quantum affine algebras of classical types (see references in \cite{h4}). For quantum toroidal algebras and more general quantum affinizations, the $q$-character of fundamental representations are totally combinatorially characterized by the Frenkel-Mukhin algorithm, and we have the $T$-system for Kirillov-Reshetikhin modules. So it should be possible to extract combinatorially explicit formulas for $q$-characters (more general than in Theorem \ref{nexp}), as well as branching rules to the underlying quantum Kac-Moody algebra.

- The structure of quantum toroidal algebras is still not fully understood for general type. In particular the role of vertical and horizontal quantum affine subalgebra and their interactions has to be described.

- It is expected, and proved for many types, that Kirillov-Reshetikhin modules (with an appropriate choice of the spectral parameter) of quantum affine algebras have a crystal basis. It would be very interesting to develop an analog theory of crystal bases for quantum toroidal algebras. As we do not have a Chevalley presentation for these algebras, this theory should rely on Drinfeld generators.

- Bezrukavnikov and Etingof defined induction and restriction functors for Cherednik algebras in \cite{be}. As suggested by P. Etingof, it would be very interesting to relate it to the fusion product of quantum toroidal algebras in a duality result analog to Proposition \ref{compt}.

- Geometric approaches to the representation theory of general quantum affinizations, in the spirit of the results of Nakajima for simply-laced quantum affinizations, have still to be developed. In particular as for the applications developed in \cite{npc, Nab}, it should be possible to interpret geometrically corresponding analogs of Kazdhan-Lusztig polynomials in order to get informations on the fusion product $*_f$. Candidates for these polynomials are constructed in \cite{h0, h1} where precise conjectures in this direction are given.


\begin{thebibliography}{99}

\bibitem[B]{bec} {\bf J. Beck}, {\it Braid group action and quantum affine algebras}, Comm. Math. Phys. {\bf 165}, no. 3, 555--568 (1994).

\bibitem[BE]{be} {\bf R. Bezrukavnikov and P. Etingof}, {\it Parabolic induction and restriction functors for rational Cherednik algebras}, {Preprint arXiv:0803.3639v2}

\bibitem[BM]{bm} {\bf W. Borho and R. MacPherson}, {\it Partial resolutions of nilpotent varieties}, {Analysis and topology on singular spaces, II, III (Luminy, 1981),  23--74, Astérisque, {\bf 101-102}, Soc. Math. France, Paris (1983)}

\bibitem[Cha1]{Chari2} {\bf V. Chari}, {\it Minimal affinizations of representations of quantum groups: the rank $2$ case}, {Publ. Res. Inst. Math. Sci.  {\bf 31},  no. 5, 873--911 (1995)}

\bibitem[Cha2]{c} {\bf V. Chari}, {\it Braid group actions and tensor
products}, {Int. Math. Res. Not. {\bf 2002}, no. 7, 357--382}

\bibitem[Che1]{che2} {\bf I. Cherednik}, {\it A unification of Knizhnik-Zamolodchikov and Dunkl operators via affine Hecke algebras}, {Invent. Math. {\bf 106} (1991),  no. 2, 411--431}

\bibitem[Che2]{che3} {\bf I. Cherednik}, {\it Double affine Hecke algebras and Macdonald's conjectures}, {Ann. of Math. (2) {\bf 141}  (1995),  no. 1, 191--216}

\bibitem[CG]{cg} {\bf V. Chari and J. Greenstein}, {\it Quantum loop modules} {Represent. Theory  {\bf 7}  (2003), 56--80}

\bibitem[CL]{cl} {\bf V. Chari and T. Le}, {\it Representations of double affine Lie algebras}, {A tribute to C. S. Seshadri (Chennai, 2002), 199--219, Trends Math., Birkhäuser, Basel, 2003}

\bibitem[CM]{cm} {\bf V. Chari and A. Moura}, {\it Characters and blocks for
finite-dimensional representations of quantum affine algebras}, {Int. Math. Res. Not. {\bf 2005}, no. 5, 257--298}

\bibitem[CP1]{cp1} {\bf V. Chari and A. Pressley}, {\it A Guide to Quantum Groups}, {Cambridge University Press, Cambridge (1994)}

\bibitem[CP2]{cp2} {\bf V. Chari and A. Pressley}, {\it Quantum affine algebras and affine Hecke algebras}, {Pacific J. Math. {\bf 174}  (1996),  no. 2, 295--326}

\bibitem[D]{Dri2}{\bf V. G. Drinfeld}, {\it A new realization of Yangians and of quantum affine algebras}, {Soviet Math. Dokl. {\bf 36}, no. 2, 212--216 (1988)}

\bibitem[E]{e} {\bf B. Enriquez}, {\it PBW and duality theorems for quantum groups and quantum current algebras}, {J. Lie Theory {\bf 13} (2003),  no. 1, 21--64}

\bibitem[EM]{em} {\bf P. Etingof and  A. Moura}, {\it Elliptic Central Characters and Blocks of Finite Dimensional Representations of Quantum Affine Algebras}, {Represent. Theory {\bf 7} (2003), 346--373}

\bibitem[FJW]{fj} {\bf I. Frenkel, N. Jing and W. Wang}, {\it Quantum vertex representations via finite groups and the McKay correspondence}, {Comm. Math. Phys. {\bf 211} (2000), no. 2, 365--393}

\bibitem[FM]{Fre2} {\bf E. Frenkel and E. Mukhin}, {\it Combinatorics of $q$-Characters of Finite-Dimensional Representations of Quantum Affine Algebras}, {Comm. Math. Phys. {\bf 216} (2001), no. 1, 23--57}

\bibitem[FR1]{frere} {\bf I. Frenkel and N. Reshetikhin}, {\it Quantum affine algebras and holonomic difference equations},  {Comm. Math. Phys.  {\bf 146},  no. 1, 1--60 (1992)}

\bibitem[FR2]{Fre} {\bf E. Frenkel and N. Reshetikhin}, {\it The $q$-Characters of Representations of Quantum Affine Algebras and Deformations of $W$-Algebras}, {Recent Developments in Quantum Affine Algebras and related topics, Cont. Math., vol. {\bf 248}, 163--205 (1999)}

\bibitem[FZ1]{fz1} {\bf S. Fomin and A. Zelevinsky}, {\it The Laurent phenomenon}, {Adv. in Appl. Math. {\bf 28} (2002),  no. 2, 119--144}

\bibitem[FZ2]{fz2} {\bf S. Fomin and A. Zelevinsky}, {\it $Y$-systems and generalized associahedra}, {Ann. of Math. (2) {\bf 158} (2003),  no. 3, 977--1018}

\bibitem[Go]{go} {\bf I. Gordon}, {\it Quiver varieties, category O for rational Cherednik algebras,
and Hecke algebras}, {Preprint arXiv:math.RT/0703150}

\bibitem[G1]{g} {\bf N. Guay}, {\it Cherednik algebras and Yangians}, {Int. Math. Res. Not. {\bf 2005}, no. 57, 3551--3593}

\bibitem[G2]{g2} {\bf N. Guay}, {\it Affine Yangians and deformed double current algebras in type A}, {Adv. Math. {\bf 211},  no. 2, 436--484 (2007)}

\bibitem[GKV]{gkv} {\bf V. Ginzburg, M. Kapranov and E. Vasserot}, {\it Langlands reciprocity for algebraic surfaces}, {Math. Res. Lett. {\bf 2} (1995), no. 2, 147--160}

\bibitem[GJ1]{gj} {\bf Y. Gao and N. Jing}, {\it Algebras over the Fock space}, {C. R. Math. Acad. Sci. Soc. R. Can. {\bf 23} (2001), no. 4, 136--140}

\bibitem[GJ2]{gj2} {\bf Y. Gao and N. Jing}, {\it $U\sb q(\widehat{gl}\sb N)$ action on $\widehat{gl}\sb N$-modules and quantum toroidal algebras}, {J. Algebra {\bf 273} (2004), no. 1, 320--343}

\bibitem[He1]{h0} {\bf D. Hernandez}, {\it Algebraic approach to q,t-characters}, {Adv. Math. {\bf 187}, no. 1, 1-52 (2004)}

\bibitem[He2]{h1} {\bf D. Hernandez}, {\it The $t$-analogs of $q$-characters at roots of unity for quantum affine algebras and beyond}, {J. Algebra {\bf 279} (2004), no. 2, 514--557}

\bibitem[He3]{h2} {\bf D. Hernandez}, {\it Representations of quantum affinizations and fusion product}, {Transform. Groups {\bf 10} (2005), no. 2, 163--200}

\bibitem[He4]{h3} {\bf D. Hernandez}, {\it Monomials of q and q,t-characters for non simply-laced quantum affinizations}, {Math. Z. {\bf 250} (2005), no. 2, 443--473}

\bibitem[He5]{hc} {\bf D. Hernandez}, {\it The Kirillov-Reshetikhin conjecture and solutions of T-systems}, {J. Reine Angew. Math. {\bf 596}  (2006), 63--87}

\bibitem[He6]{h4} {\bf D. Hernandez}, {\it On minimal affinizations of representations of quantum groups}, {Comm. Math. Phys. {\bf 277} (2007), no. 1, 221--259}

\bibitem[He7]{h5} {\bf D. Hernandez}, {\it Drinfeld coproduct, quantum fusion tensor category and applications}, {Proc. London Math. Soc. (3) {\bf 95} (2007), no. 3, 567--608}

\bibitem[He8]{h6} {\bf D. Hernandez}, {\it Smallness problem for quantum affine algebras and quiver varieties}, {Ann. Sci. \'Ecole Norm. Sup. (4) {\bf 41} (2008) , no. 2 , 271-306}

\bibitem[He9]{h8} {\bf D. Hernandez}, {\it Kirillov-Reshetikhin conjecture : the general case}, {Preprint arXiv:0704.2838}

\bibitem[Hu]{h} {\bf J. Humphreys}, {\it Reflection groups and Coxeter groups}, {Cambridge Studies in Advanced Mathematics, 29. Cambridge University Press, Cambridge (1990)}

\bibitem[HK]{hk} {\bf A. Henriques and J. Kamnitzer}, {\it The octahedron recurrence and $gl_n$ crystals}, {Adv. Math. {\bf 206},  no. 1, 211--249 (2006)}
      
\bibitem[HN]{hn} {\bf D. Hernandez and H. Nakajima}, {\it Level $0$ monomial crystal}, {Lusztig's issue of Nagoya Math. J. {\bf 184} (2006), 85--153} 

\bibitem[Jim]{j} {\bf M. Jimbo}, {\it A $q$-analogue of $U(gl(N+1))$, Hecke algebra, and the Yang-Baxter equation}, {Lett. Math. Phys. {\bf 11} (1986),  no. 3, 247--252}

\bibitem[Jin]{jin} {\bf N. Jing}, {\it Quantum Kac-Moody algebras and vertex representations}, {Lett. Math. Phys. {\bf 44} (1998),  no. 4, 261--271}

\bibitem[JKKMP]{jkkmp} {\bf M. Jimbo, R. Kedem, H. Konno, T. Miwa and J-U Petersen}, {\it Level-$0$ structure of level-$1$ $U\sb q(\widehat{{\rm sl}}\sb 2)$-modules and Macdonald polynomials} {J. Phys. A {\bf 28} (1995), no. 19, 5589--5606}

\bibitem[Kac]{kac} {\bf V. Kac}, {\it Infinite dimensional Lie algebras}, {3rd Edition, Cambridge University Press (1990)}

\bibitem[Kash1]{kasm1} {\bf M. Kashiwara}, {\it Crystal bases of modified quantized enveloping algebra},
{Duke Math. J. {\bf 73}, no. 2, 383--413 (1994)}

\bibitem[Kash2]{kas0} {\bf M. Kashiwara}, {\it On level-zero representation of quantized affine algebras},
{Duke Math. J. {\bf 112},  no. 1, 117--175 (2002)}

\bibitem[Kash3]{kas} {\bf M. Kashiwara}, {\it Realizations of crystals},
{in Combinatorial and geometric representation theory (Seoul, 2001), 133--139, Contemp. Math., {\bf 325}, Amer. Math. Soc., Providence, RI (2003)}

\bibitem[Kass]{kass} {\bf C. Kassel}, {\it Kähler differentials and coverings of complex simple Lie algebras extended over a commutative algebra}, {Proceedings of the Luminy conference on algebraic $K$-theory (Luminy, 1983). J. Pure Appl. Algebra {\bf 34},  no. 2-3, 265--275 (1984)}

\bibitem[Kn]{kn} {\bf H. Knight}, {\it Spectra of tensor products of finite-dimensional representations of Yangians}, {J. Algebra {\bf 174}  (1995),  no. 1, 187--196}

\bibitem[KTW]{ktw} {\bf A. Knutson, T. Tao and C. Woodward}, {\it A positive proof of the Littlewood-Richardson rule using the octahedron recurrence}, {Electron. J. Combin. 11 (2004), no. 1, Research Paper {\bf 61} (2004)}

\bibitem[L]{lu} {\bf G. Lusztig}, {\it Introduction to quantum groups}, {Progress in Mathematics {\bf 110} (1993) Birkhäuser Boston, Inc., Boston, MA}

\bibitem[M1]{m1} {\bf K. Miki}, {\it Toroidal braid group action and an automorphism of toroidal algebra $U_q(sl_{n+1,tor}) (n\geq 2)$}, {Lett. Math. Phys. {\bf 47} (1999), no. 4, 365--378}

\bibitem[M2]{m2} {\bf K. Miki}, {\it Representations of quantum toroidal algebra $U_q({sl}_{n+1,{tor}}) (n\geq 2)$}, {J. Math. Phys. {\bf 41} (2000), no. 10, 7079--7098}

\bibitem[M3]{m3} {\bf K. Miki}, {\it Quantum toroidal algebra $U_q({sl}_{2,{tor}})$ and $R$ matrices}, {J. Math. Phys. {\bf 42} (2001), no. 5, 2293--2308}

\bibitem[M4]{m4} {\bf K. Miki}, {\it Some quotient algebras arising from the quantum toroidal algebra $U_q({sl}_2(\CC_\gamma))$}, {Osaka J. Math. {\bf 42} (2005), no. 4, 885--929}

\bibitem[M5]{m5} {\bf K. Miki}, {\it Some quotient algebras arising from the quantum toroidal algebra $U_q({sl}_{n+1}({\CC}_\gamma)) (n\geq 2)$}, {Osaka J. Math. {\bf 43}  (2006),  no. 4, 895--922}

\bibitem[MRY]{mry} {\bf R. Moody, S. Rao, and T. Yokonuma}, {\it Toroidal Lie algebras and vertex representations}, {Geom. Dedicata {\bf 35} (1990),  no. 1-3, 283--307}

\bibitem[Nag]{nag} {\bf K. Nagao}, {\it K-theory of quiver varieties, q-Fock space and nonsymmetric Macdonald polynomials}, {Preprint arXiv:0709.1767}

\bibitem[Nak1]{Naams} {\bf H. Nakajima}, {\it Quiver varieties and finite-dimensional representations of quantum affine algebras}, {J. Amer. Math. Soc. {\bf 14} (2001), no. 1, 145--238}

\bibitem[Nak2]{npc} {\bf H. Nakajima}, {\it Geometric construction of representations of affine algebras}, {Proceedings of the International Congress of Mathematicians, Vol. I (Beijing, 2002), 423--438, Higher Ed. Press, Beijing (2002)}

\bibitem[Nak3]{Nac} {\bf H. Nakajima}, {\it $t$--analogs of $q$--characters of quantum affine algebras of type $A_n$, $D_n$}, {in Combinatorial and geometric representation theory (Seoul, 2001), 141--160, Contemp. Math., {\bf 325}, Amer. Math. Soc., Providence, RI (2003)}

\bibitem[Nak4]{Nad}{\bf H. Nakajima}, {\it $t$-analogs of $q$-characters of Kirillov-Reshetikhin modules of quantum 
affine algebras}, {Represent. Theory {\bf 7}, 259--274 (electronic) (2003)}

\bibitem[Nak5]{Nab} {\bf H. Nakajima}, {\it Quiver Varieties and $t$-Analogs of $q$-Characters of Quantum Affine Algebras}, {Ann. of Math. {\bf 160}, 1057 - 1097 (2004)} 

\bibitem[Ra]{ra} {\bf E. Rao}, {\it On representations of toroidal Lie algebras}, {Functional analysis VIII, 146--167,
Various Publ. Ser. (Aarhus), 47, Aarhus Univ., Aarhus (2004)}

\bibitem[Ro]{ro} {\bf M. Rosso}, {\it Finite-dimensional representations of the quantum analog of the enveloping algebra of a complex simple Lie algebra}, {Comm. Math. Phys. {\bf 117},  no. 4, 581--593 (1988)}

\bibitem[RR]{rr} {\bf D. Robbins and H. Rumsey}, {\it Determinants and alternating sign matrices}, {Adv. in Math. {\bf 62}, no. 2, 169--184 (1986)}

\bibitem[Sa]{Sa} {\bf Y. Saito}, {\it Quantum toroidal algebras and their vertex representations}, {Publ. Res. Inst. Math. Sci. {\bf 34} (1998), no. 2, 155--177}

\bibitem[Sc1]{Sc} {\bf O. Schiffmann}, {\it Noncommutative projective curves and quantum loop algebras}, {Duke Math. J. {\bf 121} (2004), no. 1, 113--168}

\bibitem[Sc2]{Sc2} {\bf O. Schiffmann}, {\it Canonical bases and moduli spaces of sheaves on curves}, {Invent. Math. {\bf 165} (2006), no. 3, 453--524}

\bibitem[STU]{stu} {\bf Y. Saito, K. Takemura and D. Uglov}, {\it Toroidal actions on level $1$ modules of $U_q(\hat{sl}_n)$}, {Transform. Groups {\bf 3} (1998), no. 1, 75--102}

\bibitem[TU1]{tuu} {\bf K. Takemura and D. Uglov}, {\it Level-$0$ action of $U\sb q(\widehat{ sl}_n)$ on the $q$-deformed Fock spaces}, {Comm. Math. Phys. {\bf 190} (1998),  no. 3, 549--583}

\bibitem[TU2]{tu} {\bf K. Takemura and D. Uglov}, {\it Representations of the quantum toroidal algebra on highest weight modules of the quantum affine algebra of type $gl_N$}, {Publ. Res. Inst. Math. Sci. {\bf 35} (1999), no. 3, 407--450}

\bibitem[V]{v} {\bf M. Varagnolo}, {\it Periodic modules and quantum groups}, {Transform. Groups {\bf 9}  (2004),  no. 1, 73--87}

\bibitem[VV1]{vv1} {\bf M. Varagnolo and E. Vasserot}, {\it Schur duality in the toroidal setting}, {Comm. Math. Phys. {\bf 182} (1996),  no. 2, 469--483}

\bibitem[VV2]{vv2} {\bf M. Varagnolo and E. Vasserot}, {\it Double-loop algebras and the Fock space}, {Invent. Math. {\bf 133} (1998), no. 1, 133--159}

\bibitem[VV3]{vv3} {\bf M. Varagnolo and E. Vasserot}, {\it On the $K$-theory of the cyclic quiver variety}, {Internat. Math. Res. Notices {\bf 1999}, no. 18, 1005--1028}

\end{thebibliography}
\end{document}